# A Method for Numerical Solution of Third-Kind Volterra Integral Equations Using Krall-Laguerre Polynomials

P. Jami[1], E. Hashemizadeh[*2]

*Department of Mathematics, Karaj Branch, Islamic Azad University, Karaj, Iran*

**Abstract**

The present study proposed a method for numerical solution of linear Volterra integral equations (VIEs) of the third kind, before only analytical solution methods had been discussed with reference to previous research and review of the related literature. Given that such analytical solutions are not almost always feasible, it is required to provide a numerical method for solving the mentioned equations. Accordingly, Krall-Laguerre polynomials were utilized for numerical solution of these equations. The main purpose of this method is to approximate the unknown functions through Krall-Laguerre polynomials. Moreover, an error analysis is performed on the proposed method.

## 1. Introduction

Over recent years, researchers have completed extensive scientific studies on integral equations, which can significantly contribute to modeling and analyzing a wide variety of problems in mechanics, engineering, chemistry, physics, biology, astronomy, potential theory, economics, as well as electrostatics [1- 3]. Some modeling problems have been also converted into third-kind integral equations.

For example, the following integral equation:

$$x^{\beta} f(x) = g(x) + \int_0^x (x-t)^{-\alpha} k(x,t) f(t) dt \quad , x \in [0,T] \tag{1}$$

wherein, $\alpha \in [0,1)$, $\beta \in \mathbb{R}$, $\beta > 0$, $\alpha + \beta > 0$, $g(x)$ represent continuous functions on $I$. Moreover, $k(x,t)$ is continuous on $\Delta = \{(x,t): 0 \leq t \leq x \leq T\}$ and it is in the form of $k(x,t) = x^{\alpha+\beta}k_1(x,t)$ in which $k_1 \in c(\Delta)$. This class of equations as expressed in (1) can be correspondingly found in concepts about singular integral equations with boundary value problems for mixed-kind partial differential ones. Besides, such equations have been widely applied to some major problems in the domains of neutron transport, elasticity, and scattering theory of particles [4-6]. Therefore, research studies in this area have recently received much more attention. As far as this, solutions have been investigated for existence and uniqueness of the systems of the third-kind integral equations [7-10] and researchers have provided necessary conditions to convert (1) into earlier cordial Volterra integral equations (VIEs) [11, 12]. However, the case of $\alpha + \beta > 1$ has been particularly interesting because integral operators related to (1) would not be compact and ensuring equation-solving via classical numerical methods would not be feasible if $k_1(0,0) \neq 0$. Notably, researchers [13] have presented a modified graded mesh to deal with the problem of solvability. Moreover, Gabbasov [14-16] had examined these equations through a novel direct technique and a specific correlation one. Shulaia [17, 18] had additionally studied equations on the basis of the notions of the spectral expansion solution method.

---

[1] jami.parvaneh@gmail.com
[2*] Corresponding author email: hashemizadeh@kiau.ac.ir

All these research studies had been thus conducted to find an analytical solution method for third-kind integral equations and no numerical one had been proposed. Therefore, the present study proposed a new method for numerical solution of the linear integral equations of the third kind. It should be noted that the given method was on the basis of approximating unknown functions using Krall-Laguerre polynomials.

The study was organized as follows. Section 2 deals with Krall orthogonal polynomials. Section 3 implements Krall-Laguerre polynomials on integral equations of the third kind and Section 4 discusses the estimation of error bound. In Section 5, the implementation of this method is shown on polynomials and a comparison is made between outputs obtained from this method as well as the exact solutions of these equations. At the end, the study is concluded.

## 2. Krall Orthogonal Polynomials

It should be noted that Krall orthogonal polynomials are known as subsets of polynomials with linear functional $u$ connectivity, obtained from quasi-definite functions (see [19, 20]), so $u : H \rightarrow \mathbb{N}$ ($H$ signifies the complex polynomial space with a complex coefficient) the Dirac delta function is added and $\tilde{u}$ refers to linear function.

$$\tilde{u} = u + \sum_{p=1}^{N} A_p \delta(x_p) \tag{2}$$

where, $A_p \in \mathbb{R}$, $x_1, \ldots, x_p \in \mathbb{R}$ and $\delta(x_p)$ represents Dirac delta function at point $x_p$.

For the first time in 1938, such functions were delineated by Krall [21, 22], who stated that such forms of the polynomials had emerged as essential eigen-functions of fourth-order linear differential operators with polynomial coefficients, although they were not dependent on polynomial degrees. His studies further led to Bochner's research works on development of classical orthogonal polynomials [23]. Krall also confirmed that, regardless of the classical polynomials of Hermite, Laguerre, and Jacobi, there were three additional families of orthogonal polynomials that could convince this fourth-order differential equation which were orthogonal as for measure and would not be certainly continuous according to the Lebesgue measure. For example, the Jacobi-type polynomials are orthogonal based on the weight function.

$$h(x) = (1-x)^{\alpha} + N\delta(x), \qquad N > 0,\ \alpha > -1$$

supported on $[0,1]$, Legendre-type polynomials are also orthogonal on $[-1, 1]$ according to the $h(x) = \frac{\alpha}{2} + \frac{\delta(x-1)}{2} + \frac{\delta(x+1)}{2}$, $\alpha > 0$ and Laguerre-type polynomials, which are orthogonal with regard to $h(x) = e^{-x} + N\delta(x)$, $N > 0$ on $[0, \infty)$. Such a condition resulted in further studies into orthogonal polynomials in accordance with weight functions [24, 25] in which more examples of higher-order differential equations were provided [26]. Examination of such polynomials has gained much interest in the past few years and has also led to its development [27, 28], in particular, if the starting functional $u$ is considered as one of the classical continuous linear functionals (which expanded Krall-Jacobi, Krall-Hermite, Krall-Laguerre, and Krall-Bessel polynomials) or a classical discrete one (which developed Krall-Hahn, Krall-Meixner, Krall-Kravchuk, and Krall-Charlier polynomials). The steps and how to transform Laguerre polynomials into Krall-Laguerre ones have been described in [29]. The present study used Krall-Laguerre polynomials for approximating unknown functions.

## 2.1. Krall-Laguerre Polynomials

Krall-Laguerre polynomial $K_m(x)$ of the degree $m$ is provided by [27]:

$$K_m(x) = \sum_{i=0}^{m} \frac{(-1)^i}{(i+1)!} \binom{m}{i} [i(\alpha + m + 1) + \alpha] x^i \tag{3}$$

A family of polynomials $\{K_m(x)\}_{m=0}^{\infty}$ is also orthogonal with relation to the measure $\omega$ taken by $d\omega = w(x)\,dx$, so weight function is:

$$w(x) = \frac{1}{\alpha}\delta(x) + e^{-x}H(x),$$

that, $H(x)$ is the Heaviside step function and the measure $\omega$ refers to the Laguerre weight $e^{-x}$ on $(0,\infty)$.

The first six of these polynomials are listed as follows:

$K_0(x) = 1,$

$K_1(x) = 2 - 3x,$

$K_2(x) = 2 - 7x + 2x^2,$

$K_3(x) = 2 - 12x + 7x^2 - \frac{5x^3}{6},$

$K_4(x) = 2 - 18x + 16x^2 - \frac{23x^3}{6} + \frac{x^4}{4},$

$K_5(x) = 2 - 25x + 30x^2 - \frac{65x^3}{6} + \frac{17x^4}{12} - \frac{7x^5}{120}$

The Krall-Laguerre polynomials are also depicted in Figure 1, for different $m$-values.

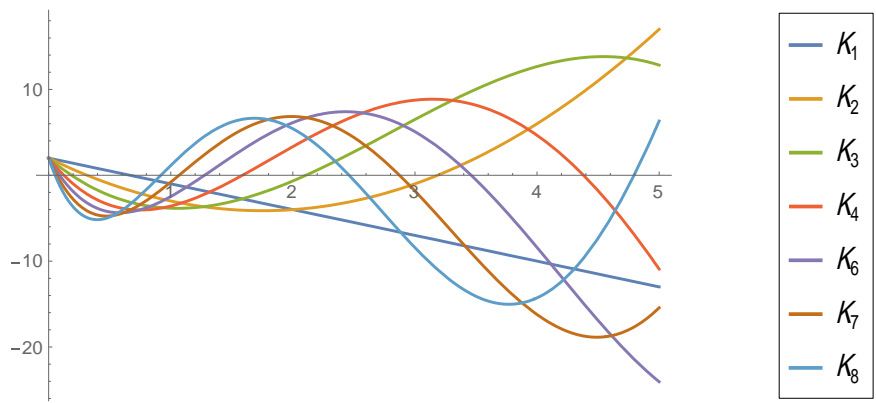

**Figure 1.** Plots of Krall-Laguerre polynomials

## 2.2. Basic Relations

In the present study, Krall-Laguerre approximation $K_m(f)$ defined to a function $f : [0,1] \to \mathbb{R}$ refers to the polynomial:

$$K_m(f)(x) = \sum_{i=0}^{m} f\left(\frac{i}{m}\right) p_{m,i}(x) \qquad (4)$$

wherein, $p_{m,i}(x)$ denotes polynomial of degree $m$:

$$p_{m,i}(x) = \frac{(-1)^i}{(i+1)!} \binom{m}{i} [i(\alpha+m+1)+\alpha] x^i \qquad i = 0, \ldots, m. \qquad (5)$$

It should be noted that the error bound of the Krall-Laguerre polynomial is examined in the next section.

## 3. Implementation of Krall-Laguerre Method for Third-Kind Integral Equation

Here, the third-kind integral equation would be presented as follows:

$$x^{\beta} f(x) = g(x) + \int_0^x (x-t)^{-\alpha} k(x,t) f(t) dt, \quad x \in [0,T] \qquad (6)$$

To numerically solve these kinds of integral equations, the unknown $f$ is approximated by relation (4). Accordingly, this equation would be provided:

$$\sum_{i=0}^{m} f\left(\frac{i}{m}\right) \left( \frac{(-1)^i}{(i+1)!} [i(\alpha+m+1)+\alpha] \binom{m}{i} \left( \left(x^{\beta+i}\right) - \int_0^x (x-t)^{-\alpha} k(x,t) t^i dt \right) \right) = g(x)$$

To find $f\left(\frac{i}{m}\right)$, $i = 0, \ldots, m$, the equation would be converted into a linear system of equations through replacing $x$ with $x_j = j/m + \varepsilon$, $j = 0, \ldots, m$, $j = 0, \ldots, m$, and $x_m = 1 - \varepsilon$, so $\varepsilon$ would be arbitrary little. It should be noted that each distinct value in $[0, 1]$ would be chosen as $x_j$, $j = 0, \ldots, m$ except for the singular values of in the integral equation in this study. Therefore, the singularity could not be skipped. After that, the equation could be written as:

$$BX = Y \tag{7}$$

where,

$$B = \left[\frac{(-1)^i}{(i+1)!}[i(\alpha+m+1)+\alpha]\binom{m}{i}\left((x_j^{\beta+i}) - \int_0^x (x_j-t)^{-\alpha} k(x_j,t) t^i dt\right)\right], \quad i,j=0,1,\ldots,m,$$

$$X = \left[f\left(\frac{i}{m}\right)\right]^t, \quad i=0,\ldots,m, \tag{8}$$

$$Y = \left[g(x_j)\right]^t, \quad j=0,\ldots,m,$$

thus, the integral would be computed, which is found in the $B$'s formula numerically. Here, $f\left(\frac{i}{m}\right)$, $i=0,\ldots,m$, by $f_m\left(\frac{i}{m}\right)$, $i=0,\ldots,m$, could be shown, which are here considered as solutions in the nodes. Replacing them in Equation (4), $k_m(f_m)(x_i)$, $i=0,\ldots,m$, is obtained that is the solution for integral Equation (6).

## 4. Error Analysis

To get started, the following definition is considered:

The space $W_r^p$ is the weighted Sobolev space of order $p$.

$$W_r^P = \{f \in L_r^2 \mid \|f\|_{W_s^p} < +\infty\}$$

That norm is described by:

$$\|f\|_{W_s^p} = \left(\int_{\mathbb{R}_+} f^2(x) w(x) dx\right)^{\frac{1}{2}}$$

**Theorem 1.** [28] Assume $f(x) \in W_s^P$, $1 < p < +\infty$ and $f_m(x) = \sum_{i=0}^m a_i k_i(x)$ is the best approximation polynomial of $f(x)$ in $L_2-$ norm, then, for each $0 \leq s$:

$$\|f - f_m\|_{W_r^p} \leq C \frac{\log m}{m^s} \|f\|_{W_s^p}, \quad p \geq 1. \tag{9}$$

For purpose of convenience in calculations, Equation (1) would be transformed into Fredholm integral equation of the second type, so:

$$F(x) = x^{\beta} f(x),$$

Then,

$$F(x) = g(x) + \int_0^x (x-t)^{-\alpha} k(x,t) f(t) dt, \quad x \in [0, \infty)$$

Now, the unknowns are approximated by Krall-Laguerre polynomials.

$$K_m(F(x)) = \tilde{g}(x) + \int_0^x (x-t)^{-\alpha} k(x,t) K_m(f(t)) dt$$

In general:

$$g(x) = F(x) + \int_0^x (x-t)^{-\alpha} k(x,t) f(t) dt,$$

so,

$$\tilde{g}(x) = K_m(F(x)) - \int_0^x (x-t)^{-\alpha} k(x,t) K_m(f(t)) dt$$

In Equation (7), if $X$ is replaced with its approximation using Krall-Laguerre polynomial as shown by $X_n$, then:

$$BX_n = \tilde{Y} \tag{10}$$

If $B$ is invertible, thenceforth:

$$\| X - X_n \| \leq B^{-1} \| Y - \tilde{Y} \|$$

in which, $X_n$ is the solution obtained from the system.

So, there is a need to compute the error bound of $\| Y - \tilde{Y} \|$

$$\sup | Y - \tilde{Y} | = \sup | g(x_i) - \tilde{g}(x_i) |$$

$$\sup \left| (F(x_i) - K_m(F(x_i)) - \int_0^x (x_i - t)^{-\alpha} k(x_i, t)((f(t) - K_m(f(t))) dt \right|$$

$$\leq \sup | F(x_i) - K_m(F(x_i)) | + \sup \left| \int_0^x (x_i - t)^{-\alpha} k(x_i, t)(f(t) - K_m(f(t))) dt \right|$$

$$\leq \sup\left| (F(x_i) - K_m(F(x_i))\right| + \sup\left| (x_i - t)^{-\alpha} k(x_i,t)(f(t) - K_m(f(t))\right| = \blacklozenge$$

assume that $x_i \in (0, \infty)$ and $\mu = \sup|(x-t)^{-\alpha} k(x,t)|$

Using Theorem 1, then:

$$\blacklozenge \leq C \frac{\log m}{m^s} \|f\|_{W_r^p} + C\mu \frac{\log m}{m^s} \|f\|_{W_r^p}$$

$$= (1+\mu) C \frac{\log m}{m^s} \|f\|_{W_r^p}$$

$$\Rightarrow \|X - X_n\| \leq \|B^{-1}\|(1+\mu) C \frac{\log m}{m^s} \|f\|_{W_r^p}$$

It should be noted that the extra condition was used in the lemma below. Thus, a bound for $\|B^{-1}\|$ and a condition number of $B$ are shown.

**Lemma 1.** It is assumed that $\|B - I\| = C_0 < 1$,

$\|.\|$ refers to the highest norm of the rows and $I$ represents identity matrix of order $n+1$.

Therefore, $\|B^{-1}\| \leq \dfrac{1}{1-C_0}$ and:

$$\text{cond}(B) \leq \frac{1+C_1}{1+C_0},$$

where, $C_1$ is $\max_i \left(\int_0^{x_i} |(x_i - t)^{-\alpha} k(x_i, t)| dt\right)$.

**Proof.** A bound is determined for $\|B\|$. Matrix $B$ would be known from the relation, thus:

$$\|B\| = \max_j \sum_{i=0}^{n} \left| \frac{(-1)^i}{(i+1)!}[i(\alpha+m+1)+\alpha]\binom{m}{i}\left((x_j^{\beta+i}) - \int_0^{x_j}(x_j - t)^{-\alpha} k(x_j, t) t^i dt\right)\right|$$

$$= \max_j \left|1 - \int_0^{x_j}(x_j - t)^{-\alpha} k(x_j - t) dt\right| \leq (1 + C_1).$$

Accordingly, a bound is required for $\|B^{-1}\|$.

$$\|D\| = \|B - I\| = C_0 < 1.$$

Then,

$$\|B^{-1}\| = \|(I+D)^{-1}\| \leq \frac{1}{1-\|D\|},$$

As a result,

$$\text{cond}(B) = \|B\|\|B^{-1}\| \leq \frac{1+C_1}{1+C_0},$$

and then the proof is completed.

## 5. Illustrative Examples

Now, two different instances are shown, denoting that the given method can be accurate, applicable, and effective. According to [29], presenting both examples before, the efficiency of the method proposed in this study was evaluated. So, it is assumed that $T = 1$ and $h = \frac{1}{m}$ is expressed.

To reflect on method error, notations would be introduced:

$$e_m = \max_{0 \leq i \leq n} |f(t_i) - f_m(t_i)|,$$

$$p_m = \log_2\left(\frac{e_m}{e_{2m}}\right)$$

Wherever $t_i = ih$, $f(t)$ is the accurate solution, $f_m(t)$ is the approximate solution to the suggested method. Therefore, calculations were fulfilled on a personal computer with a Core-i7 processor, 2.40 GHz frequency, and 8 GB Memory, and then the codes were expressed in Mathematica.

## 6. Example of Third-Kind Integral Equation

**Example 1**. With regard to the first example, Equation (1) would be considered with $\alpha = \frac{2}{3}, \beta = \frac{2}{3}, k(t,x) = \frac{\sqrt{3}}{3\pi} x^{1/3}$ that can provide an equation of Abel type:

$$x^{\frac{2}{3}} f(x) = g(x) + \int_0^x \frac{\sqrt{3}}{3\pi} t^{\frac{1}{3}} (x-t)^{-\frac{2}{3}} f(t) \, dt, \quad t \in [0,1]$$

where,

$$g(x) = x^{\frac{47}{12}} \left( 1 - \frac{\Gamma\left(\frac{1}{3}\right)\Gamma\left(\frac{55}{12}\right)}{\pi\sqrt{3}\Gamma\left(\frac{59}{12}\right)} \right)$$

Put differently, the exact solution of the equation would be $f(x) = x^{\frac{13}{4}}$.

First, the above equation was solved via various m-values. Then, numerical outputs were listed in Table 1 reporting the highest error, order of convergence, and outputs of the collocation method [11]. Notably, numerical outputs demonstrated that the proposed method had a convergence of 3.06. Whereas, in the collocation method described in [11], these examples had an order of convergence by 2.96.

**Table 1.** Numerical outputs for example 1

| Methods | m | $e_m$ | $P_m$ |
|---|---|---|---|
| Present method (i.e. Krall-Laguerre method) | 3 | $4.12 \times 10^{-6}$ | ------ |
|  | 6 | $1.02 \times 10^{-4}$ | 2.01 |
|  | 12 | $1.22 \times 10^{-5}$ | 3.06 |
|  | 24 | $1.32 \times 10^{-6}$ | 3.21 |
|  | 48 | $1.37 \times 10^{-7}$ | 3.27 |
|  | 96 | $1.45 \times 10^{-8}$ | 3.24 |
| Method developed in [11] |  |  |  |
| N-Points ($n = 3$) | 256 | $4.58 \times 10^{-7}$ | 2.96 |
| Radau II ($n = 3$) | 256 | $5.13 \times 10^{-9}$ | 3.28 |

**Example 2.** Consider now Equation (1) with $\beta=1$, $\alpha=0$, $k(t,x)=\frac{1}{2}$, that would be applied for modeling a number of heat conduction problems with the mixed-kind boundary conditions:

$$x f(x) = \frac{6}{7} x^3 \sqrt{x} + \int_0^x \frac{1}{2} f(t) \, dt, \qquad t \in [0,1]$$

The above equation would have the accurate solution of $f(x) = x^{\frac{5}{2}}$.

The present method with its various m-values was utilized in this equation. Table 2 also reports the outputs. The highest error and the order of convergence have been correspondingly indicated. Then, a comparison has been made and outputs attained from the collocation method [11] are outlined in Table 2. Accordingly, the order convergence for the proposed method has been estimated to be 4.49 with regard to the one provided in [11]. It should be noted that the order of convergence for the above-mentioned example would be 1.99.

**Table 2.** Computed errors of $o\ e_n$ for example 2

| Method | $m$ | $e_m$ | $P_m$ |
|---|---|---|---|
| Present method (Krall-Laguerre) | 3 | $7.55 \times 10^{-4}$ | ------ |
|  | 6 | $3.37 \times 10^{-5}$ | 4.49 |
|  | 12 | $3.50 \times 10^{-6}$ | 4.49 |
|  | 24 | $9.13 \times 10^{-8}$ | 5.26 |
| Method [11] | | | |
| Chebyshev (n=2) | $N = 256$ | $1.30 \times 10^{-5}$ | 1.99 |

## 7. Conclusion and Further Research

The present study proposed a numerical method on the basis of Krall-Laguerre polynomials for solving third-kind VIEs. This method was developed for reducing the problem to the systems of algebraic equations that are not readily solvable. Therefore, two numerical examples were provided to confirm the applicability and precision of the method. Thus, the numerical outputs showed that the order of convergence of this method was appropriate.

There are a few methods for solving integral equations of the third kind; therefore, they can be solved by other numerical methods in future work and their results can be compared with the ones in the present study.